\makeatletter
\def\input@path{{style_files}}
\makeatother
\documentclass{article}

\usepackage{style_files/arxiv}
\usepackage{amssymb}
\usepackage{amsmath}
\usepackage{empheq}
\usepackage{amscd}
\usepackage{graphicx}
\usepackage{graphics}
\usepackage[noadjust]{cite}
\usepackage{amsthm}
\usepackage{caption}
\usepackage{multirow}
\usepackage{array}
\usepackage{rotating}

\usepackage{indentfirst}
\usepackage{tikz}
\usepackage{calc}
\usepackage{epsfig}
\usepackage{natbib}
\usepackage[makeroom]{cancel}
\usepackage{multicol}
\usepackage{csquotes}
\usepackage{algorithm}
\usepackage{algorithmic}
\usepackage{enumitem}
\usepackage{hyperref}
\usepackage{url}
\hypersetup{colorlinks,linkcolor={blue},citecolor={blue},urlcolor={blue}}

\usepackage{style_files/timet}
\usepackage{graphicx}
\usepackage{multirow}
\usepackage{multicol}
\usepackage{lipsum}
\usepackage{epsfig}
\usepackage{amsmath}
\usepackage{amsfonts}
\usepackage{amssymb}
\usepackage{color}
\numberwithin{equation}{section}
\usepackage{caption}
\usepackage{float}
\usepackage{subcaption}
\usepackage{graphics}
\usepackage{xcolor,graphicx}
\usepackage{csquotes}
\usepackage{mathtools}
\usepackage{style_files/optidef}

\newcommand{\vertiii}[1]{{\left\vert\kern-0.25ex\left\vert\kern-0.25ex\left\vert #1 
    \right\vert\kern-0.25ex\right\vert\kern-0.25ex\right\vert}}

\usepackage{hyperref}       

\renewcommand{\P}[0]{\bold{P}}                        

\DeclareMathOperator*{\minimize}{minimize}

\begin{document}
\title{Multipliers waveform inversion}

\author{
  \href{https://orcid.org/0000-0002-9879-2944}{\includegraphics[scale=0.06]{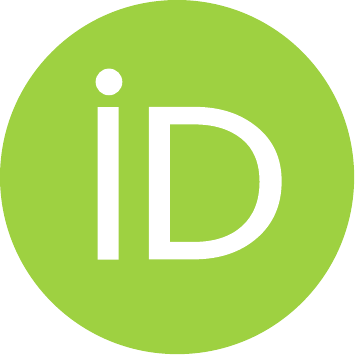}\hspace{1mm}Ali Gholami} \\
  Institute of Geophysics, University of Tehran, Tehran, Iran.
  \texttt{agholami@ut.ac.ir} \\ 
\And
\href{http://orcid.org/0000-0003-1805-1132}{\includegraphics[scale=0.06]{style_files/orcid.pdf}\hspace{1mm}Hossein S. Aghamiry} \\
  University Cote d'Azur - CNRS - IRD - OCA, Geoazur, Valbonne, France. 
  \texttt{aghamiry@geoazur.unice.fr}
  \And
\href{http://orcid.org/0000-0002-4981-4967}{\includegraphics[scale=0.06]{style_files/orcid.pdf}\hspace{1mm}St\'ephane Operto} \\ 
  University Cote d'Azur - CNRS - IRD - OCA, Geoazur, Valbonne, France. 
  \texttt{operto@geoazur.unice.fr}
  }

\renewcommand{\shorttitle}{Multipliers waveform inversion, Gholami et al.}

\maketitle

\begin{abstract}
The full waveform inversion (FWI) addresses the computation and characterization of subsurface model parameters by matching predicted data to observed seismograms in the frame of nonlinear optimization. 
We formulate FWI as a nonlinearly constrained optimization problem, for which a regularization term is minimized subject to the nonlinear data matching constraint.
Unlike FWI which is based on the penalty function, the method of multipliers solves the resulting optimization problems by using the augmented Lagrangian function; and leads to a two step recursive algorithm. The primal step requires solving an unconstrained minimization problem like the traditional FWI with a difference that the data are replaced by the Lagrange multipliers.
The dual step involves an update of the Lagrange multipliers. The overall performance of the algorithm is improved considering that this multiplier method does not require exact solve of these primal-dual subproblems. In fact, convergence is attained when only one step of a gradient-based  method is taken on both subproblems. 
The proposed algorithm greatly improves the overall performance of FWI such as convergence from inaccurate starting models and robustness with respect to determination of the step length.
Furthermore, it can be performed by the existing FWI engines with minimal change. We only have to replace the observed data at each iteration with the multipliers, thus all the nice properties of the traditional FWI algorithms are kept.
Numerical experiments confirm that the multipliers waveform inversion can converge to a solution of the inverse problem in the absence of low-frequency data from an inaccurate initial model even with a constant step size.
\end{abstract}
\graphicspath{{"./figures/"}}
The full waveform inversion (FWI) has been widely accepted as an accurate method for the computation and characterization of subsurface model parameters by matching predicted data to observed seismograms \citep{Tarantola_1988_TBI,Pica_1990_NIS,Pratt_1998_GNF,Virieux_2009_OFW}. 
The FWI involves the solution of a nonlinearly constrained  minimization problem, involving an objective function \citep[such as a measure of the roughness of the subsurface parameters like first or second derivatives, the total variation (TV), or a combination of them][]{Aghamiry_2019_CRO} and two constraints: a linear data matching constraint and a nonlinear wave-equation constraint. The nonlinearity of the wave-equation arises due to the coupling between the model parameters and the wavefield. The equation is linear with respect to the wavefield and nonlinear with respect to the model parameters. 
Taking advantage of this special structure, traditional FWI methods eliminate the wavefield (the linear variable) obtaining a somewhat more complicated forward map that involves only the nonlinear model parameters, a technique referred to as the variable projection method \citep{Golub_2013_VPM}.
The resulting reduced nonlinear problem is then recast as an equivalent quadratic penalty form, suitable to be solved by unconstrained optimization tools \citep{Nocedal_2006_NO}. 
Due to the non-linearity of the FWI a natural approach for solving it relay on first or second order gradient-based numerical methods such as the preconditioned steepest descent, nonlinear conjugate gradients, or Newton-type methods \citep{Metivier_2017_TRU}. The penalty formulation, however, has some well-known disadvantages such as sensitivity to the accuracy of the initial model, sensitivity to the penalty parameter choice, and slow convergence rate \citep{Bertsekas_1976_OPM}. 
The former item causes the FWI algorithm to trap in a local minimum when the misfit function is nonconvex and the initial model is not in the basin including the global minimum. This happens as long as the events in the predicted data (calculated from the initial model) are away by more than half a cycle from the corresponding events in the observed data, a phenomenon called cycle skipping \citep{Virieux_2009_OFW}.
Thus the most straightforward way for solving FWI is to use initial models of sufficient accuracy. Such an initial model may be obtained from, e.g., traveltime tomography, slope tomography, or a combination of them \citep{Tavakoli_2017_STE,Sambolian_2019_BIM}. 
However, building an accurate initial model is itself a challenging ill-posed problem.

During the past decade, considerable attention has been given to devising methods for mitigating cycle skipping and thus increasing the robustness of the FWI to initial models \citep[see the review paper by][]{Virieux_2009_OFW}. 
A simple classical strategy is to perform the inversion using a frequency continuation strategy (multiscale inversion) because the forward modelling map is more
linear at lower frequencies \citep{Bunks_1995_MSW}. While effective, this approach fails when the data lack low-frequency content, which is the case in most practical applications. Following this, \citet{Li_2016_FWI} extrapolated the missing low frequencies of the data to
improve the accuracy of the initial model. \citet{Treister_2017_FWI} augmented FWI with traveltime tomography to directly inject the low-frequency information of the latter into the inversion. 
Some researches replaced the usual quadratic function with alternative functions to measure the misfit (the distance between predicted data and observed data) in a hope to decrease the nonconvexity of the problem. For example, correlation or deconvolution based misfit functions \citep{VanLeeuwen_2010_CMC,Luo_2011_DBO,Warner_2016_AWI} and optimal transport misfit function \citep{Yang_2017_AOT,Metivier_2019_GOT}.

Some approaches, referred to as extended FWI, consider adding extra (physical or nonphysical) degrees of freedom to the search space \citep{Symes_2008_MVA}. The added parameters allow more flexibility to fit the data at early (or even first) iterations to a sufficient accuracy. Then, through iterations, these algorithms try to remove the effect of the extension by adjusting the model (while keeping the misfit at the desired level) such that it turn to the original reduced problem at the convergence point.  
 \citet{VanLeeuwen_2013_MLM} do this by relaxing the wave-equation constraint by using a penalty formulation of the unreduced FWI. The extra parameter introduced is the so called data-assimilated wavefield \citep{Aghamiry_2020_AED} and the try is to approach this wavefield to the original reduced wavefield by forcing the wave-equation constraint. \citet{Huang_2018_VSE} introduced nonphysical sources and tried to collapse them into the physical source through iterations.  
\citet{Metivier_2020_ARE} proposed relaxing the receiver location at first then forced to converge to the corresponding physical receiver location at the convergence point.

In this paper, we use the method of multipliers, also known as augmented Lagrangian method, to solve the reduced FWI problem.
The method was proposed and named in 1969 by \citet{Hestenes_1969_MAG} and was also proposed independently in the same year by  \citet{Powell_1969_NLC} for effective resolution of nonlinearly constrained optimization problems.
As a combination of the Lagrangian and penalty methods, the method of multipliers have global convergence guarantees under relatively weak assumptions, converge considerably faster than the pure penalty methods, and are relatively intensive to the penalty parameter choice \citep[see][ for a complete comparison between penalty and multiplier methods]{Bertsekas_1976_OPM}.
Furthermore, multiplier methods are particularly useful for solving large-scale nonlinearly constrained optimization problems because they can be implemented matrix-free \citep{Miele_1972_OTM,Tapia_1977_DMM,Nocedal_2006_NO}. 
\citet{Aghamiry_2019_IWR} used a special variant of the method of multipliers, known as alternating-direction method of multipliers \citep{Boyd_2011_DOS}, to solve unreduced FWI and gained improvements over the wavefield reconstruction inversion of  \citet{VanLeeuwen_2013_MLM} which is based on the penalty formulation of the same problem.
The main drawback is a need to compute the data-assimilated wavefields at each iteration and also storage of the multipliers vector associated with the wave-equation constraint. These may limit somewhat its applicability for time-domain formulations \citep{Wang_2016_FIR,Rizzuti_2019_ADF,Aghamiry_2020_AED,Gholami_2020_EFW}.

In this paper, we apply the method of multipliers directly to the reduced FWI where the wavefield is completely eliminated through variable projection. This leads to a minimax problem with minimization over the model parameters (the primal variable)  and maximization over Lagrange multipliers (the dual variable). We remind that the dual variable in this case is of the size of data. 
The algorithm solves this minimax problem by using alternating optimization. At each iteration, the objective function is first minimized with respect to the model parameters with multipliers held fixed, then the multiplier are updated for maximizing the dual objective. 
The first subproblem, as given in the original formulation given by \citet{Hestenes_1969_MAG} and \citet{Powell_1969_NLC}, requires us to solve a classical FWI in which the data are replaced with the multipliers vector.
However, following a diagonalized variant of the method  \citep{Miele_1972_OTM,Tapia_1977_DMM} we design an optimal algorithm by performing only a single iteration of a gradient-based method to partially minimize the objective and hence increase the efficiency while maintaining the desirable properties of the algorithm. 
Also, we simply use a first order gradient ascent method to update the multipliers, as proposed in the original papers by \citet{Hestenes_1969_MAG,Powell_1969_NLC}.
Although both the primal and dual subproblems can be performed by a Newton-type update for rapid convergence of the algorithm, we show that a simple preconditioned steepest descent update of the model and a gradient ascent update of the multipliers lead to a good progress towards a minimizer even for complicated velocity models like the 2004 BP salt model. However, performing this steps by Newton-type family can improve the result.

\section{Theory}
The full waveform inversion (FWI) with a general regularization term can be formulated as the following PDE constrained optimization problem \citep{Aghamiry_2019_CRO}: 
\begin{mini} 
{{\bold{m}\in \mathcal{B},\bold{u}}}{\mathcal{R}(\bold{m})}
{\label{main_optim}}{}
\addConstraint {\bold{A(m)u}}{=\bold{b}^*}
\addConstraint {\bold{Pu}}{=\bold{d}^*},
\end{mini}
where $\mathcal{R}$ is the (convex) regularization function, which incorporates the prior information about the subsurface parameters $\bold{m}$,  
$\mathcal{B}=\{\bold{m}\vert \bold{m}_{min} \leq \bold{m} \leq \bold{m}_{max}\}$ in which $\bold{m}_{min}$ and $\bold{m}_{max}$ denote the lower bound and upper bound vectors for the desired model satisfying $\bold{m}_{min}\leq\bold{m}_{max}$, 
$\bold{u}\in {\mathbb{R}}^{N\times 1}$ is the wavefield, $\bold{A(m)} \in {\mathbb{R}}^{N \times N}$ is the PDE operator, $\bold{b}^* \in {\mathbb{R}}^{N\times 1}$ is the source term, $\bold{d}^* \in {\mathbb{R}}^{M \times 1}$ is the recorded seismic data, and the observation operator $\bold{P} \in {\mathbb{R}}^{M\times N}$ (with $M\ll N$) samples $\bold{u}$ at receivers.
More specifically, $N=N_x\times N_z \times N_t$ and $M=N_t\times N_r$, where $N_x$ and $N_z$ denote the number of grid points sampling the subsurface model $\bold{m}$ in the horizontal and vertical dimensions, respectively, $N_t$ is the number of time samples, and $N_r$ is the number of receivers.
In this paper, we limit ourselves to the scalar acoustic wave equation with constant density equal to 1, namely, $\bold{A(m)} = \bold{m} \circ \partial_{tt} - \boldsymbol{\nabla}^2$, where $\bold{m}$ denotes the squares slowness, $\partial_{tt}$ the second time derivative and $\boldsymbol{\nabla}^2$ the Laplace operator.
Furthermore, we derive our formulation for single-source problems. However, the extension to multiple sources follows easily by summation over the sources.

The FWI optimization problem given by equation \ref{main_optim} is nonconvex. The nonconvexity in FWI arises due to the presence of the bilinear term  $\bold{m} \circ \partial_{tt}\bold{u}$ in the wave equation constraint. Different methods have been proposed to solve FWI in either time or frequency domains \citep{Tarantola_1988_TBI,Pica_1990_NIS,Pratt_1998_GNF,Virieux_2009_OFW}. 
Each method consider solving the problem in either the original unreduced-space or reduced-space \citep{Golub_2013_VPM}. The reduced-space optimization problem is obtained by eliminating the state variable $\bold{u}$ from the problem, which drastically reduces the dimensionality of the problem compared to the unreduced-space:
\begin{mini} 
{\bold{m}\in \mathcal{B}}{\mathcal{R}(\bold{m})}
{\label{main_optim_red}}{}
\addConstraint {\bold{S}(\bold{m})\bold{b}^*}{=\bold{d}^*}
\end{mini}
where $\bold{S(m)}=\bold{PA(m)^{-1}}$. This dimensionality reduction is gained at the cost of treating with more complicated forward operator $\bold{S(m)}$ and loosing the biconvexity of the optimization. However, the reduced problem can converge in fewer iterations than if the same minimization
algorithm is used for the unreduced problem and the cost per iteration for the reduced problem can be lower compared to that for the full problem \citep{Golub_2013_VPM}.

Both unreduced and reduced problems can be solved by different optimization algorithms. 
In most cases, the algorithms seeks the solution by replacing the original constrained problem by a sequence of unconstrained subproblems \citep{Nocedal_2006_NO}.
The simplest and widely used method, the quadratic penalty method, replaces the constraints by penalty terms in the objective function, where each penalty term is a scale (the penalty parameter) of the Euclidean distance of the constraint violation. 
For example, the classical FWI is obtained by a penalty formulation of equation \ref{main_optim_red}:
\begin{equation} \label{penalty}
\minimize_{\bold{m}\in \mathcal{B}}~ \mathcal{R}(\bold{m})+\frac{\mu}{2}\|\bold{S}(\bold{m})\bold{b}^* -\bold{d}^*\|_2^2,
\end{equation}
where $\mu>0$ is the penalty parameter. Equation \ref{penalty} is the most widely used formulation for FWI, either with a regularization term \citep[e.g.,][]{Gao_2019_AEW,Aghamiry_2020_FWI} or simply focusing on the resolution of the feasibility problem by setting $\mathcal{R}=0$ \citep[e.g.,][]{Gauthier_1986_TDN,Pica_1990_NIS,Pratt_1998_GNF,Virieux_2009_OFW,Metivier_2017_TRU}. We remind that, in practice, some kind of damping is applied to the gradient for determination of a suitable step direction which may be explained by including a regularization term in the objective function even though it is not mentioned explicitly.

The penalty formulation of the unreduced FWI, equation \ref{main_optim}, is
\begin{equation} \label{wri}
\minimize_{\bold{m}\in \mathcal{B},\bold{u}}~ \mathcal{R}(\bold{m})+\frac{\mu}{2}\|\bold{Pu} -\bold{d}^*\|_2^2 + \frac{\gamma}{2}\|\bold{A}(\bold{m})\bold{u} -\bold{b}^*\|_2^2,
\end{equation}
with two penalty parameters $\mu$ and $\gamma$. \citet{Esser_2018_TVR} consider solving this problem with $\mathcal{R}$ as the total-variation (TV) regularization while 
\citet{VanLeeuwen_2013_MLM} consider solving the feasibility problem ($\mathcal{R}=0$).

Despite their appeal, penalty methods have some important disadvantages; the main difficulty lies in choosing appropriate values of the penalty parameter. Their convergence rate can also be considerably slow \citep{Bertsekas_1976_OPM}.
The method of multipliers or the augmented Lagrangian method  \citep{Hestenes_1969_MAG} was proposed to overcome the shortcomings of the penalty method and now it is one of the most effective methods for solving the optimization problems  with nonlinear constraints \citep{Hestenes_1969_MAG,Powell_1969_NLC,Nocedal_2006_NO}.  
Three main advantages of the multiplier method are: 
1) convergence from inaccurate/remote starting points, 
2) efficiency for solving large-scale nonlinear optimization problems because they can be implemented matrix-free, and 
3) robustness with respect to the penalty parameter choice. The augmented Lagrangian formulation of the unreduced FWI, equation \ref{main_optim}, is
\begin{equation} \label{ir-wri}
\min_{\bold{m}\in \mathcal{B},\bold{u}}\max_{\boldsymbol{\lambda},\boldsymbol{\nu}}~ \mathcal{R}(\bold{m})-\langle\boldsymbol{\lambda},\bold{Pu} -\bold{d}^*\rangle-\langle\boldsymbol{\nu},\bold{A}(\bold{m})\bold{u} -\bold{b}^*\rangle+\frac{\mu}{2}\|\bold{Pu} -\bold{d}^*\|_2^2 + \frac{\gamma}{2}\|\bold{A}(\bold{m})\bold{u} -\bold{b}^*\|_2^2,
\end{equation}
where $\boldsymbol{\lambda}$ and $\boldsymbol{\nu}$ are Lagrange multiplier vectors associated with the data and the wave-equation constraints and $\langle\cdot,\cdot\rangle$ denotes inner product.
Recently, the authors have used the method of multipliers to solve equation \ref{ir-wri} either without the regularization term \citep{Aghamiry_2019_IWR} or with different forms of the regularization function \citep{Aghamiry_2019_IBC,Aghamiry_2019_CRO,Aghamiry_2020_FWI}.
Even though the alternating direction strategy (alternating direction method of multipliers (ADMM), \citealp*{Boyd_2011_DOS}) allowed breaking the primal optimization over the full-space into smaller subproblems over each parameter class separately, the algorithm face two main challenges from the computation and memory point of views:
1) reconstruction of the data-assimilated wavefield, as required at each iteration of the optimization algorithm, is time consuming and 
2) storage of the Lagrange multipliers associated to the wave-equation constraint is memory demanding because it is of the size of the wavefield. 
These challenges prevent the algorithm to be fully applied in the time domain applications \citep{Aghamiry_2020_AED} or even in the frequency-domain for large scale  applications where direct methods can not be used to solve the data-assimilated system.
Several attempts have been made to improve these limitations. In a recent work \citep{Gholami_2020_EFW}, the authors, using a data space formulation, formulated the wavefield reconstruction problem in the data space to reduce the computational burden of the wavefield reconstruction. Furthermore, they projected the Lagrange multipliers associated with the wave equation in the data space at the expense of additional adjoint simulations. 
Considering promising results obtained by the method of multipliers in solving the unreduced FWI, equation \ref{main_optim}, in the next section, we apply this powerful method to solve the reduced FWI, equation \ref{main_optim_red}. To the best of our knowledge, this has not been considered before.

\section{FWI by the method of multipliers}
In this section, we propose an iterative algorithm for solving the FWI problem in equation \ref{main_optim_red}. 
Ideally, we would like to compute a global minimizer of the problem when the problem is feasible, even though, it is well acknowledged that estimation of even a local minimizer is computationally difficult. In the case that the problem is infeasible, we would like to find a solution that at least satisfies the first-order optimality condition of the nonlinear feasibility problem, i.e., where the gradient of the misfit term vanish. 

The augmented Lagrangian formulation of the reduced FWI in equation \ref{main_optim_red} reads
\begin{equation}
\min_{\bold{m}\in \mathcal{B}} \max_{\boldsymbol{\lambda}} ~\mathcal{R}(\bold{m})- \langle\boldsymbol{\lambda},\bold{S}(\bold{m})\bold{b}^* -\bold{d}^*\rangle+\frac{\mu}{2}\|\bold{S}(\bold{m})\bold{b}^* -\bold{d}^*\|_2^2,
\end{equation}
where $\boldsymbol{\lambda}$ is the Lagrangian multiplier vector associated with the constraint.
This augmented Lagrangian function can be interpreted as the Lagrangian function of equation \ref{main_optim_red} when the objective is replaced by the penalty function, equation \ref{penalty}. 
The method of multipliers solves this min-max problem iteratively. At each iteration of this algorithm, the objective function
is first minimized with respect to $\bold{m}$ with $\boldsymbol{\lambda}$ held fixed,  then the multiplier $\boldsymbol{\lambda}$ is updated for maximizing the objective. The simplest role of multiplier update, as suggested in the original papers by \citet{Hestenes_1969_MAG} and \citet{Powell_1969_NLC} leads to the following iteration:
\begin{align} 
\bold{m}_{k+1} &=\arg\min_{\bold{m}\in \mathcal{B}} ~\mathcal{R}(\bold{m})- \langle\boldsymbol{\lambda}_k,\bold{S}(\bold{m})\bold{b}^* -\bold{d}^*\rangle+\frac{\mu}{2}\|\bold{S}(\bold{m})\bold{b}^* -\bold{d}^*\|_2^2, \label{MMP}\\
\boldsymbol{\lambda}_{k+1} &=\boldsymbol{\lambda}_k - \mu(\bold{S}(\bold{m}_{k+1})\bold{b}^*-\bold{d}^*), \label{MMD}
\end{align}
where $k$ is the iteration count and the simplest choose for the initial multiplier is $\boldsymbol{\lambda}_0=\bold{0}$. Motivation for this update role can be gained by differentiating the augmented Lagrangian objective with respect to the model parameters:
\begin{equation} \label{ALgrad}
\frac{\partial}{\partial \bold{m}}\mathcal{R}(\bold{m}) - \bold{J(m)}^t[\boldsymbol{\lambda}_k - \mu(\bold{S}(\bold{m})\bold{b}^*-\bold{d}^*)],
\end{equation}
where
\begin{equation} \label{J}
\bold{J(m)}=\frac{\partial \bold{S(m)}\bold{b}^*}{\partial \bold{m}}=-\bold{P}\bold{A(m)}^{-1}\frac{\partial \bold{A(m)}}{\partial \bold{m}}\bold{A(m)}^{-1}\bold{b}^*,
\end{equation}
is the Jacobian matrix and superscript $t$ denotes matrix transposition.
We can see that the right most term in the brackets, equation \ref{ALgrad}, can be considered as the Lagrangian multiplier vector that partially satisfies
the first order optimality condition of the original constrained problem, equation \ref{main_optim_red}.
This shows that a reasonable choice for the multiplier vector at next iteration should be that is given in equation \ref{MMD}.
Also, if we use a variable projection method to recast the augmented Lagrangian function in equivalent dual form which involves only the multiplier vector then the gradient of the dual function will be $-(\bold{S}(\bold{m}(\boldsymbol{\lambda}))\bold{b}^*-\bold{d}^*)$ and thus the optimum multiplier satisfies $\bold{S}(\bold{m}(\boldsymbol{\lambda}^*))\bold{b}^*=\bold{d}^*$ \citep{Luenberger_2010_LNP}. 
Accordingly, the update role in equation \ref{MMD} can be interpreted as a gradient ascent step (with the step length  $\mu$) applied to the dual problem to locally maximize it with respect to $\boldsymbol{\lambda}$.
 Considering this, one may apply Newton-like methods on this dual problem to improve the convergence at the expense of increased computational complexities \citep{Tapia_1977_DMM}.

The primal subproblem in equation \ref{MMP} can be further simplified by completing the squares and ignoring the constant term involving $\boldsymbol{\lambda}_k$, which gives (by a change of variable $\bold{d}_k\leftarrow\boldsymbol{\lambda}_k/\mu$)
\begin{align} 
\bold{m}_{k+1} &=\arg\min_{\bold{m}\in \mathcal{B}} ~\mathcal{R}(\bold{m})+\frac{\mu}{2}\|\bold{S}(\bold{m})\bold{b}^* -\bold{d}_k\|_2^2, \label{MWI1}\\
\bold{d}_{k+1} &=\bold{d}_k + \bold{d}^* - \bold{S}(\bold{m}_{k+1})\bold{b}^*. \label{MWI2}
\end{align}
This scaled form iteration begins by setting $\bold{d}_0=\bold{d}^*$.
Interestingly, the model subproblem is exactly the classical FWI in equation \ref{penalty} in which the data are replaced with the multipliers, hence multipliers waveform inversion (MWI).
The MWI as given by equations \ref{MWI1}-\ref{MWI2} involves the resolution of a sequence of FWI problems for inverting the multipliers waveform, each involving an estimation of new multipliers. 
There has been substantial activity in the area of FWI by gradient-based numerical methods (e.g. the preconditioned steepest descent, the l-BFGS quasi-Newton, Gauss-Newton, or full-Newton algorithms). The results of previous research on FWI thus can be employed for MWI.
At first glance, the MWI may be extremely costly and inefficient because it requires multiple FWI resolutions. 
However, the optimal formulation of MWI will require only one step of the mentioned iterative methods for approximate minimization of equation \ref{MWI1} before the multipliers (and possibly the penalty parameter) are updated. The convergence of the algorithm with such inexact resolution of the primal subproblem was considered by \citet{Miele_1972_OTM}, and the algorithm was termed as diagonalized variant of the original method \citep{Tapia_1977_DMM}. The digonalization approach has been widely used for solving linear and nonlinear problems, e.g., \citet{Goldstein_2009_SBM,Gholami_2017_CNA,Aghamiry_2019_IBC,Dokht_2020_AIA}.
In the following subsections we provide the detailed algorithm for performing MWI. 
\subsection{Update of model}
Since the model subproblem, equation \ref{MWI1}, is nonlinear, the standard approach for solving it is iterative linearization. 
Considering the characteristics of the regularization function $\mathcal{R}$, different methods can be used for solving this subproblem. A recent paper \citep{Aghamiry_2020_FWI} discusses an algorithm based on proximal Newton method for a general form (black box) regularization function. In solving such nonlinear problem a large class of algorithms approximate the nonlinear misfit function $E(\bold{m})=\|\bold{S}(\bold{m})\bold{b}^* -\bold{d}_k\|_2^2$ by a simple quadratic function near a reference model  $\bold{m}_k$:
\begin{equation} \label{quad}
E(\bold{m})
\approx E_k
+ \langle\nabla E_k,\bold{m-m}_k\rangle + \frac{1}{2}(\bold{m-m}_k)^t\bold{H}_k(\bold{m-m}_k)
\end{equation}
where $E_k, \nabla E_k$, and $\bold{H}_k$ are, respectively, the value, gradient, and (approximate) Hessian of $E$, all evaluated at $\bold{m}_k$. 
The gradient of $E$ is 
\begin{equation} \label{grad1}
\nabla E = \bold{J}^t(\bold{S}(\bold{m})\bold{b}^* -\bold{d}_k),
\end{equation}
where $\bold{J}$ is the Jacobian matrix defined in equation \ref{J}. There are many strategies for choosing $\bold{H}$. If we choose it to be the exact Hessian of the misfit term, $\bold{H}=\nabla^2 E$, then we get the Newton approximation of the misfit function and hence Newton MWI. 
Generally, the strategies for choosing Hessian approximations in Newton-type FWI methods can also be adapted to MWI.
Three most common choices are as follows:
\begin{enumerate}
\item[1)] A scale of the identity matrix, the most simple choice:
\begin{equation} \label{H0}
\bold{H}_k=\alpha_k\bold{I},
\end{equation}
which specializes to the gradient descent when $\mathcal{R}=0$. $\alpha_k$ serves as the step length which can be determined at each iteration by a line search method.
\item[2)] The diagonal pseudo-Hessian matrix \citep{Shin_2001_IAP}: 
\begin{equation} \label{hess1}
 \bold{H}_k=\alpha_k \bold{L}_k^t\bold{L}_k,
\end{equation} 
where
\begin{equation}
\bold{L}_k=\frac{\partial \bold{A(m)}}{\partial \bold{m}}\bold{A}(\bold{m}_k)^{-1}\bold{b}^*.
\end{equation}
In this case, the method specializes to a preconditioned gradient descent when $\mathcal{R}=0$. 
This diagonal Hessian is simple to implement and accounts for geometrical spreading. 
\item[3)] Damped Gauss-Newton Hessian:
\begin{equation} \label{hess2}
 \bold{H}_k=\bold{J}_k^t\bold{J}_k + \epsilon \bold{I},
\end{equation}
which specializes to the damped Gauss-Newton or Levenberg-Marquardt method when $\mathcal{R}=0$ \citep{Levenberg_1944_MSC}. 
Parameter $\epsilon$ is a positive scalar used for stabilizing the Hessian inverse. 
This choice of Hessian is more favourable compared with the previous ones because it accounts for the curvature of misfit function more accurately and is free of tuning step length. Even though some people also include a step length with the Newton-type methods, we do not induce it here because even with a fine-tuned step length the solution of the subproblem is still approximate. We remind that, due to great defect correction action of the multipliers, MWI algorithm is convergent with approximate model updates. For practical applications, the construction and application of the Gauss-Newton Hessian would seem to be prohibitively large and computationally expensive \citep{Metivier_2017_TRU}. 
\end{enumerate}

\subsubsection{Efffcient computation of Gauss-Newton MWI}
In order to simplify the computation of the model update with the quadratic form in equation \ref{quad} while using the Gauss-Newton Hessian in equation \ref{hess2}, we choose a modified gradient vector and the diagonal pseudo-Hessian in our approximation such that the global minimizer of the new quadratic function coincides with that of the Gauss-Newton approximation.
More specifically, we choose the following quadratic function:
\begin{equation} \label{quad_GND}
E(\bold{m})
\approx E_k
+ \langle\bold{J}_k^t\bold{Q}_k^{-1}(\bold{S}_k\bold{b}^*-\bold{d}_k),\bold{m-m}_k\rangle + \frac{1}{2}\|\bold{L}_k(\bold{m-m}_k)\|_2^2
\end{equation}
where $\bold{Q}_k=\bold{S}_k\bold{S}_k^t + \epsilon \bold{I}$ is the data domain part of the Gauss-Newton Hessian. 
For surface acquisitions where the number of receivers is smaller than the number of model parameters (matrix  $\bold{S}$ is under-determined) this approach will reduce the size of the Hessian matrix dramatically.

Some advantages of the new formulation are: 1) The remaining Hessian which applies in the image space is diagonal. This brings computational advantages for minimization of the linearized model subproblem. In this case, a number of algorithms can be used to efficiently solve it for a general regularization function \citep[see][ for more details]{Aghamiry_2019_CRO,Aghamiry_2020_FWI}. 2) The data domain Hessian, $\bold{S}\bold{S}^t + \epsilon \bold{I}$, is source independent and we need to construct it only ones and use it for multiple sources. This is different from the data domain Hessian proposed by \citet[][ Their equation 40]{Pratt_1998_GNF} which is source dependent.

\section{Numerical examples}
In this section, we demonstrate the performance of MWI in comparison with FWI on synthetic models including unscaled 2004 BP salt model while using remote starting models and in the absence of low-frequency data. In all tests, we use the preconditioned (applied with the pseudo-Hessian matrix in equation \ref{hess1}) gradient descent for model updating. A fixed step length is used in all tests. 

\subsection{Transmission example}
\begin{figure}
\center
\includegraphics[width=1\columnwidth]{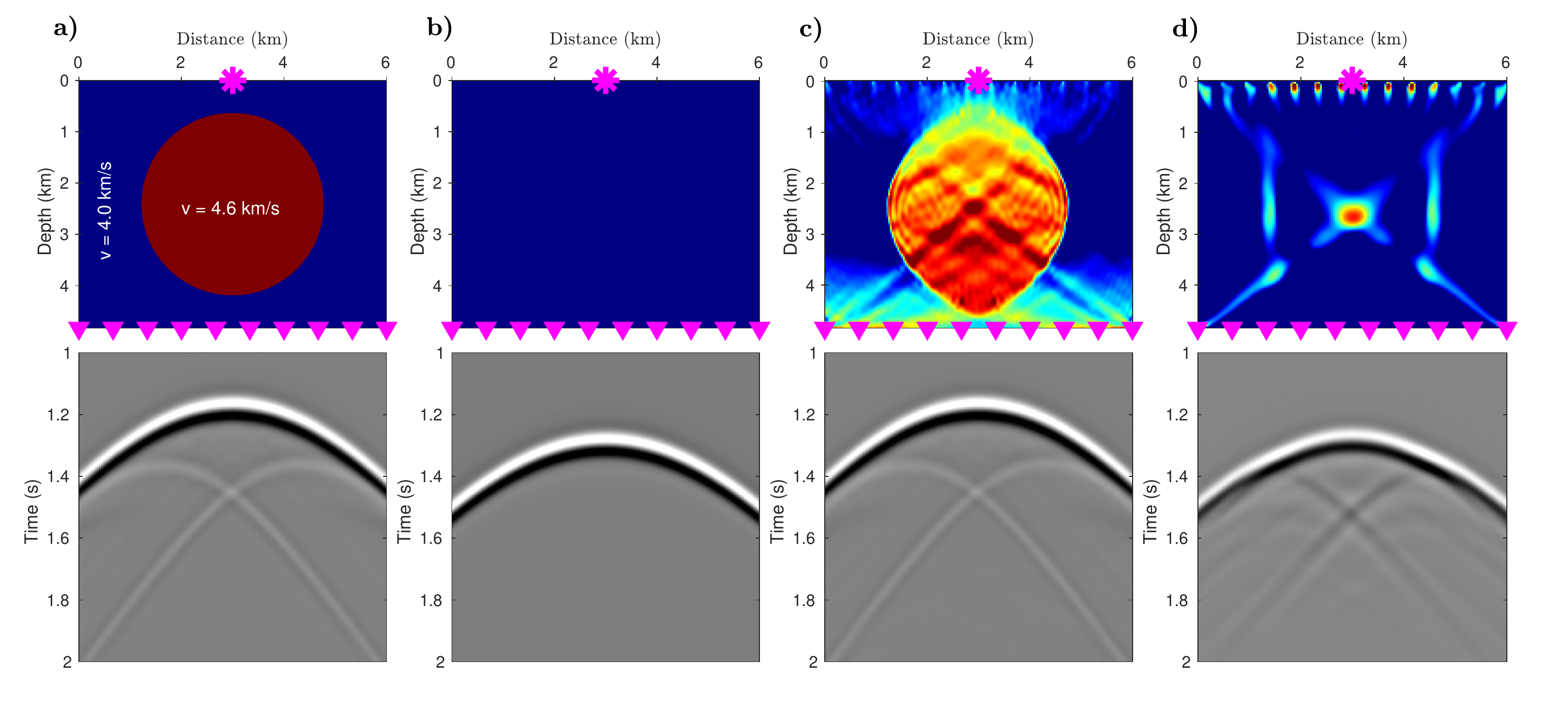}
\caption{Camembert transmission example. (a) true velocity. (b) initial velocity. (c-d) velocities inferred from (c) MWI and (d) FWI. Bottom row shows the associated shot gathers recorded by receivers (triangles) for the sources (stars).}
\label{Camembert}
\end{figure}

Here we test the performance of MWI in comparison with FWI by using the ``Camembert" transmission example. 
The subsurface model contains a circular anomaly of velocity 4.6 km/s embedded in a homogeneous background of velocity 4.0 km/s (Figure \ref{Camembert}a). 
The dimensions of the model are 4.8 km in distance and 6 km in depth, and the grid spacing is 35.5 m.
The acquisition consists of 14 equally spaced sources in the top side of the model and 170 equally spaced receivers at the bottom side (shown by triangles). The source signature is a 10~Hz Ricker wavelet. The shot gather associated with the source at distance 3 km is shown bellow the velocity model.
We start the inversion from the homogeneous background model (Figure \ref{Camembert}b) to meet the nonlinear regime of the FWI. The cycle skipping of the first arrivals can be seen clearly by comparing the associated shot gathers shown bellow each velocity model.
We perform one cycle of the inversion with the 10 Hz Ricker wavelet (one cycle means that the full bandwidth is inverted in one go).
The results obtained using MWI and FWI are presented in Figures \ref{Camembert}c and \ref{Camembert}d, respectively. 
We can observe from the velocity models and the sample shots in Figure \ref{Camembert} that MWI tackled the cycle skipping and recovered the true velocity model successfully while FWI trapped in a local minimum for picking a meaningless velocity.

\subsection{Reflection example}
\begin{figure}
\center
\includegraphics[width=1\columnwidth]{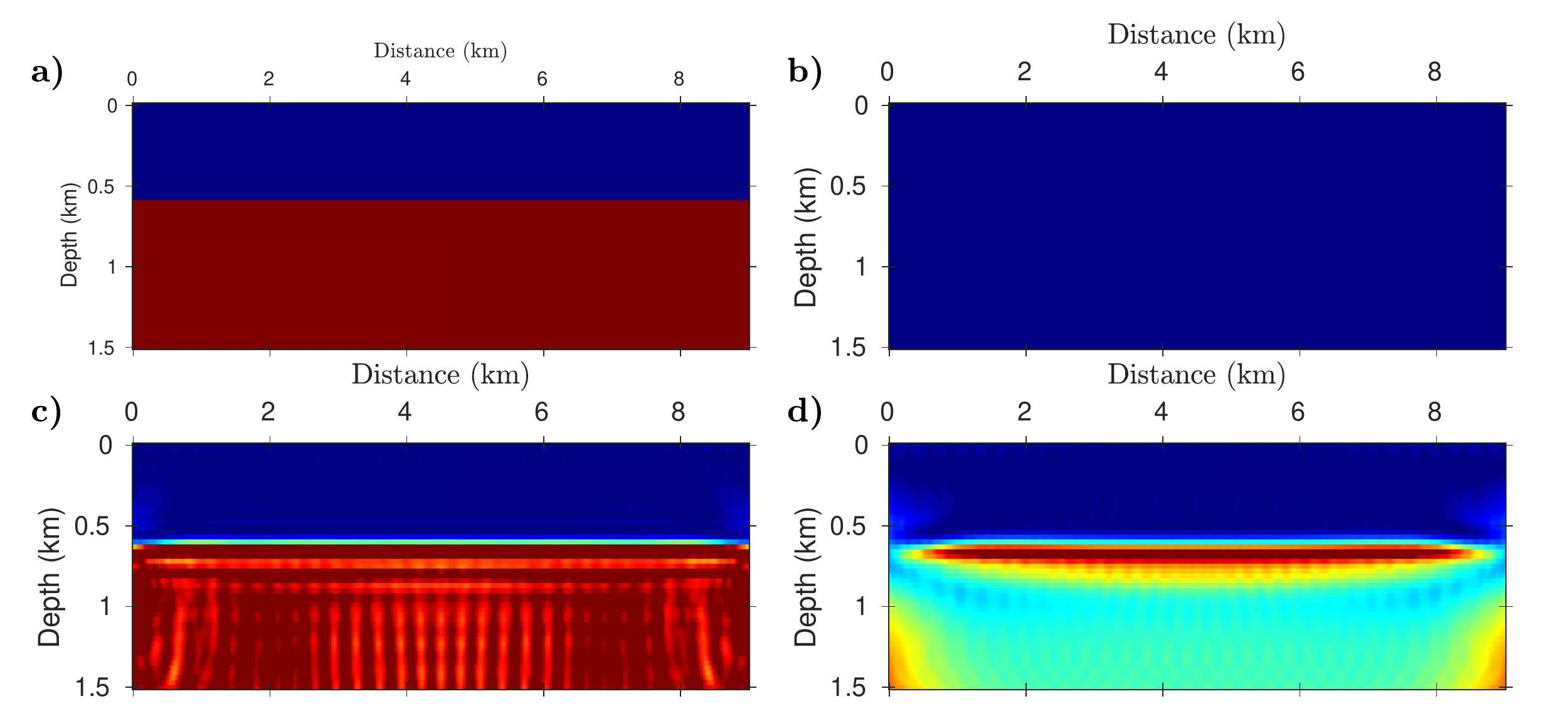}
\caption{Reflection example. (a) true velocity. (b) initial velocity. (c-d) velocities inferred from (c) MWI and (d) FWI. Surface acquisition consist of 31 sources and 310 receivers positioned uniformly on the top side of the model.}
\label{Reflection_vels}
\end{figure}

We continued by testing the performance of MWI in comparison with FWI by reconstructing a velocity model which contains two layers of homogenous constant velocity from surface data.
The dimensions of the model are 9 km in distance and 1.5 km in depth, and the grid spacing is 30 m.
The velocity model contains two layer with thickness 0.6 and 0.9 km with constant velocity 2.0 and 4.0 km/s, respectively (Figure \ref{Reflection_vels}a).
The fixed spread acquisition consists of 31 equally spaced sources and 310 equally spaced receivers in the top side of the model. The source signature is a 5~Hz Ricker wavelet.
The total recording time is 1 second. PML boundary conditions are included all around the model so that reflections from the interface are the only information from the second layer which are recorded by receivers and used for the reconstruction.
We start the inversion from the homogeneous initial velocity 2.0 km/s (Figure \ref{Reflection_vels}b) and perform one cycle of the inversion with the 5 Hz Ricker wavelet.
The results obtained using MWI and FWI (after 400 iterations while no regularization and bound constraints are implemented) are presented in Figures \ref{Reflection_vels}c and \ref{Reflection_vels}d, respectively.
The MWI was able to successfully reconstruct the velocity  of second layer with correct interface but the FWI results mainly provided the the information about the central part of the layer interface while the majority of the information about the layer velocity is missed.

\subsection{2004 BP salt model}
\begin{figure}
\center
\includegraphics[width=1\columnwidth]{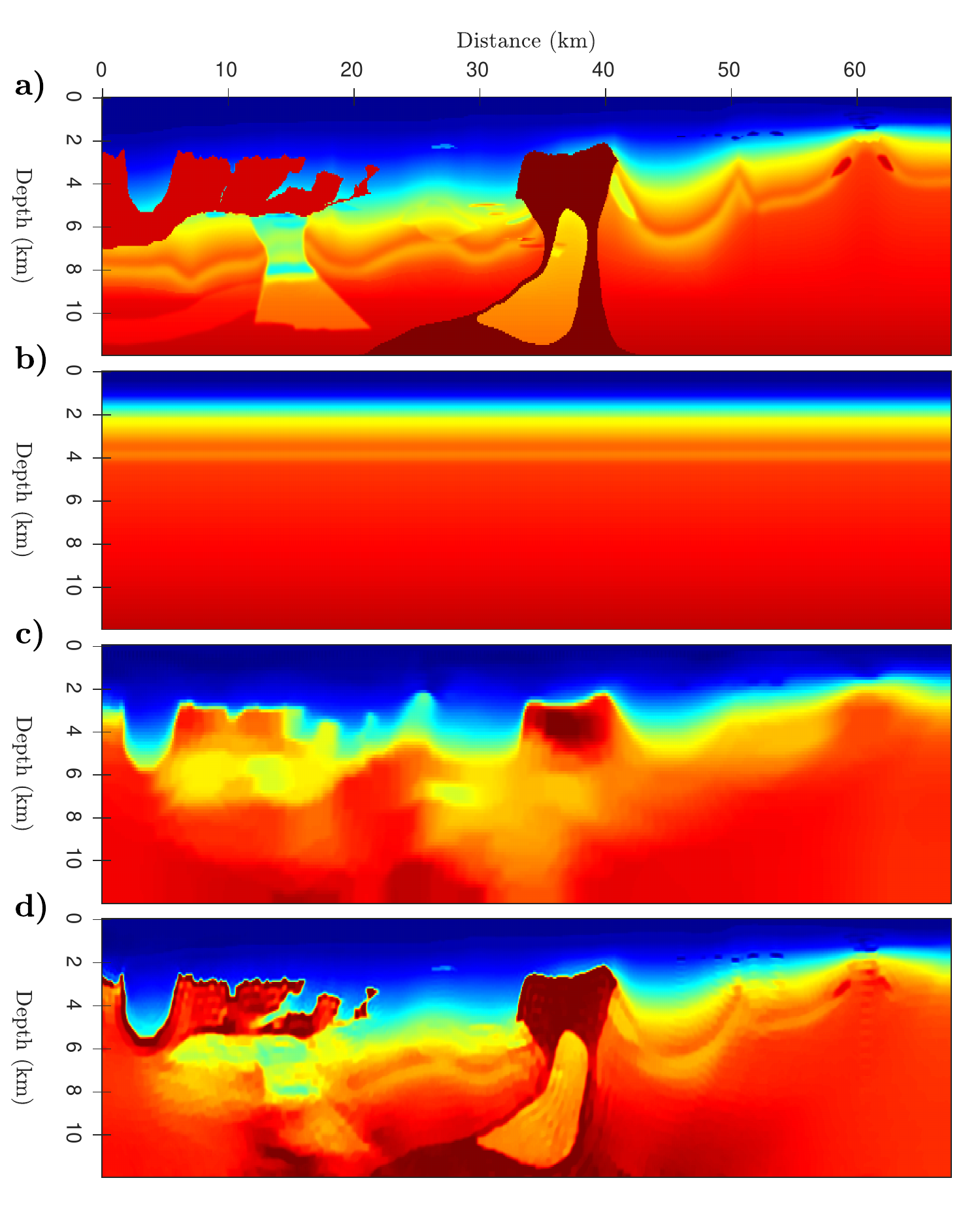}
\caption{The 2004 BP salt example: (a) true velocity. (b) 1D initial velocity.
(c) MWI inverted macro model with grid size 150 m by successive inversion of frequencies 1.0 Hz, 1.5 Hz, 2.0 Hz, and 2.5 Hz (bound constrained TV regularization included). (d) final MWI result with grid size 75 m, started with (c) by successive mono-frequency inversion of 1-4 Hz.}
\label{BP2004}
\end{figure}

Finally, we assess the proposed MWI method against the original unscaled 2004 BP salt model when a crude starting model is used (Figure \ref{BP2004}a).
The subsurface model is 67.5 km wide and 12 km deep. It contains two salt bodies with large velocity contrast and steeply dipping flanks. Also, there exist several low-velocity anomalies in the shallow part of the model. 
The fixed spread acquisition consists of equally-spaced 225 sources and equally-spaced 450 receivers on top side of the model. 
The source signature is a Ricker wavelet with a 3 Hz dominant frequency. 
We discretized the subsurface model to 160 $\times$ 900 grid points with a 75 m grid interval and used a very erroneous 1D initial model (Figure \ref{BP2004}b).   
Starting from this initial point the highest starting frequency in the data to avoid cycle-skipping on this model is 0.3 Hz \citep{Sun_2020_EFWI}. 
Therefore, the starting frequency 1 Hz leads to significant cycle-skipping that drive FWI to a local minimum. 
We used seven frequencies from 1 Hz to 4 Hz with frequency interval of 0.5 Hz in the inversion. 
To reduce the computation burden, we further downsample the model to 80 $\times$ 450 grid points with a grid interval of 150 m.
Then we used two cycles of single-frequency inversion moving from the 1 Hz frequency to the 2.5 Hz frequency according to a classic frequency continuation
strategy. However, to partially prevent inverse crime, we inverted the same data generated from the fine gridded model.
We also included total variation regularization and bound constraint (with the minimum and maximum velocity of the true model) at this stage. The final macro model generated at this stage is shown in Figure \ref{BP2004}c. We can observe that the algorithm successfully overcame the cycle skipping and generated the low-wavenumber structures with true shapes and kinematics. Then, we upsampled Figure \ref{BP2004}c by a factor of two and performed other two cycles of single-frequency inversion, this time, moving from 1 to 4 Hz frequencies on fine gridded model.  Figure \ref{BP2004}d shows the final resulting model.
We can observe that the  recovered model by MWI for this challenging problem is reasonably good.

%
%

\section{Conclusions}
Full waveform inversion (FWI) was formulated as nonlinearly constrained optimization problem which requires minimization of a regularization term subject to a nonlinear constraint that ensures matching predicted data to observed seismograms. The classical reduced space FWI was obtained by the penalty formulation of this problem. However, it converges to a local minimum as long as the predicted data are away by more that half a cycle from the observed seismograms. 
In this paper, we solved the original problem by using the multiplier method (augmented Lagrangian method), termed as multipliers waveform inversion (MWI).
MWI replaces the observed data in FWI by the Lagrange multiplies which are the running sum of the data residuals in iterations based on a simple gradient ascent update applied to the dual objective. The Lagrange multiplies progressively correct the adjoint source in the gradient with respect to the model parameters of the misfit function. Numerical examples were presented showing that the MWI can converge to a solution of the inverse problem in the absence of low-frequency data from an inaccurate initial model even with a constant step size.

\bibliographystyle{style_files/gji}

\newcommand{\SortNoop}[1]{}

\end{document}